\documentclass[11pt]{amsart}
\usepackage{amssymb,mathrsfs}
\title{A Note on A Fundamental Domain for Siegel-Jacobi Space}

\author{Jae-Hyun Yang}
\address{Department of Mathematics, Inha University,
Incheon 402-751, Korea}
\email{jhyang@inha.ac.kr }

\begin{document}

\newtheorem{theorem}{Theorem}[section]
\newtheorem{lemma}{Lemma}[section]
\newtheorem{proposition}{Proposition}[section]
\newtheorem{remark}{Remark}[section]
\newtheorem{definition}{Definition}[section]

\renewcommand{\theequation}{\thesection.\arabic{equation}}
\renewcommand{\thetheorem}{\thesection.\arabic{theorem}}
\renewcommand{\thelemma}{\thesection.\arabic{lemma}}
\newcommand{\BR}{\mathbb R}
\newcommand{\BQ}{\mathbb Q}
\newcommand{\bn}{\bf n}
\def\charf {\mbox{{\text 1}\kern-.24em {\text l}}}
\newcommand{\BC}{\mathbb C}
\newcommand{\BZ}{\mathbb Z}

\thanks{
\indent This work was supported by INHA UNIVERSITY Research Grant
(INHA-31619)}


\begin{abstract}
{In this paper, we study a fundamental domain for the
Siegel-Jacobi space $Sp(g,{\mathbb Z})\ltimes H_{\mathbb
Z}^{(g,h)}\backslash {\mathbb H}_g\times {\mathbb C}^{(h,g)}$.}
\end{abstract}
\maketitle


\newcommand\tr{\triangleright}
\newcommand\al{\alpha}
\newcommand\be{\beta}
\newcommand\g{\gamma}
\newcommand\gh{\Cal G^J}
\newcommand\G{\Gamma}
\newcommand\de{\delta}
\newcommand\e{\epsilon}
\newcommand\z{\zeta}
\newcommand\vth{\vartheta}
\newcommand\vp{\varphi}
\newcommand\om{\omega}
\newcommand\p{\pi}
\newcommand\la{\lambda}
\newcommand\lb{\lbrace}
\newcommand\lk{\lbrack}
\newcommand\rb{\rbrace}
\newcommand\rk{\rbrack}
\newcommand\s{\sigma}
\newcommand\w{\wedge}
\newcommand\fgj{{\frak g}^J}
\newcommand\lrt{\longrightarrow}
\newcommand\lmt{\longmapsto}
\newcommand\lmk{(\lambda,\mu,\kappa)}
\newcommand\Om{\Omega}
\newcommand\ka{\kappa}
\newcommand\ba{\backslash}
\newcommand\ph{\phi}
\newcommand\M{{\Cal M}}
\newcommand\bA{\bold A}
\newcommand\bH{\bold H}

\newcommand\Hom{\text{Hom}}
\newcommand\cP{\Cal P}
\newcommand\cH{\Cal H}

\newcommand\pa{\partial}

\newcommand\pis{\pi i \sigma}
\newcommand\sd{\,\,{\vartriangleright}\kern -1.0ex{<}\,}
\newcommand\wt{\widetilde}
\newcommand\fg{\frak g}
\newcommand\fk{\frak k}
\newcommand\fp{\frak p}
\newcommand\fs{\frak s}
\newcommand\fh{\frak h}
\newcommand\Cal{\mathcal}

\newcommand\fn{{\frak n}}
\newcommand\fa{{\frak a}}
\newcommand\fm{{\frak m}}
\newcommand\fq{{\frak q}}
\newcommand\CP{{\mathcal P}_g}
\newcommand\Hgh{{\mathbb H}_g \times {\mathbb C}^{(h,g)}}
\newcommand\BD{\mathbb D}
\newcommand\BH{\mathbb H}
\newcommand\CCF{{\mathcal F}_g}
\newcommand\CM{{\mathcal M}}
\newcommand\Ggh{\Gamma_{g,h}}
\newcommand\Chg{{\mathbb C}^{(h,g)}}
\newcommand\Yd{{{\partial}\over {\partial Y}}}
\newcommand\Vd{{{\partial}\over {\partial V}}}

\newcommand\Ys{Y^{\ast}}
\newcommand\Vs{V^{\ast}}
\newcommand\LO{L_{\Omega}}
\newcommand\fac{{\frak a}_{\mathbb C}^{\ast}}

%
%
\begin{section}{{\bf Introduction}}
\setcounter{equation}{0}
For a given fixed positive integer $g$, we let
$${\mathbb H}_g=\,\{\,\Omega\in \BC^{(g,g)}\,|\ \Om=\,^t\Om,\ \ \ \text{Im}\,\Om>0\,\}$$
be the Siegel upper half plane of degree $g$ and let
$$Sp(g,\BR)=\{ M\in \BR^{(2g,2g)}\ \vert \ ^t\!MJ_gM= J_g\ \}$$
be the symplectic group of degree $g$, where $F^{(k,l)}$ denotes
the set of all $k\times l$ matrices with entries in a commutative
ring $F$ for two positive integers $k$ and $l$, $^t\!M$ denotes
the transpose matrix of a matrix $M$ and
$$J_g=\begin{pmatrix} 0&I_g\\
                   -I_g&0\end{pmatrix}.$$
$Sp(g,\BR)$ acts on $\BH_g$ transitively by
\begin{equation}
M\cdot\Om=(A\Om+B)(C\Om+D)^{-1},
\end{equation} where $M=\begin{pmatrix} A&B\\
C&D\end{pmatrix}\in Sp(g,\BR)$ and $\Om\in \BH_g.$ Let $\G_g$ be
the Siegel modular group of degree $g$. C. L. Siegel\,[8] found a
fundamental domain ${\mathcal F}_g$ for $\G_g\ba\BH_g$ and
calculated the volume of $\CCF.$ We also refer to [2],\,[4],\,[10]
for some details on $\CCF.$ \vskip 0.1cm For two positive integers
$g$ and $h$, we consider the Heisenberg group
$$H_{\BR}^{(g,h)}=\{\,(\la,\mu;\ka)\,|\ \la,\mu\in \BR^{(h,g)},\ \kappa\in\BR^{(h,h)},\
\ka+\mu\,^t\la\ \text{symmetric}\ \}$$ endowed with the following
multiplication law
$$(\la,\mu;\ka)\circ (\la',\mu';\ka')=(\la+\la',\mu+\mu';\ka+\ka'+\la\,^t\mu'-
\mu\,^t\la').$$ We define the semidirect product of $Sp(g,\BR)$
and $H_{\BR}^{(g,h)}$
$$G^J=Sp(g,\BR)\ltimes H_{\BR}^{(g,h)}$$
endowed with the following multiplication law
$$
(M,(\lambda,\mu;\kappa))\cdot(M',(\lambda',\mu';\kappa')) =\,
(MM',(\tilde{\lambda}+\lambda',\tilde{\mu}+ \mu';
\kappa+\kappa'+\tilde{\lambda}\,^t\!\mu'
-\tilde{\mu}\,^t\!\lambda'))$$ with $M,M'\in Sp(g,\BR),
(\lambda,\mu;\kappa),\,(\lambda',\mu';\kappa') \in
H_{\BR}^{(g,h)}$ and
$(\tilde{\lambda},\tilde{\mu})=(\lambda,\mu)M'$. Then $G^J$ acts
on $\BH_g\times \BC^{(h,g)}$ transitively by
\begin{equation}
(M,(\lambda,\mu;\kappa))\cdot (\Om,Z)=(M\cdot\Om,(Z+\lambda
\Om+\mu)
(C\Omega+D)^{-1}), \end{equation} where $M=\begin{pmatrix} A&B\\
C&D\end{pmatrix} \in Sp(g,\BR),\ (\lambda,\mu; \kappa)\in
H_{\BR}^{(g,h)}$ and $(\Om,Z)\in \BH_g\times \BC^{(h,g)}.$ We note
that the Jacobi group $G^J$ is {\it not} a reductive Lie group and
also that the space ${\mathbb H}_g\times \BC^{(h,g)}$ is not a
symmetric space. We refer to [11]-[14] and [16] about automorphic
forms on $G^J$ and topics related to the content of this paper.
From now on, we write $\BH_{g,h}:=\BH_g\times \BC^{(h,g)}.$ \vskip
0.1cm We let $$\G_{g,h}:=\G_g\ltimes H_{\BZ}^{(g,h)}$$ be the
discrete subgroup of $G^J$, where
$$H_{\BZ}^{(g,h)}=\{\,(\la,\mu;\ka)\in H_{\BR}^{(g,h)}\,|\ \la,\mu\in \BZ^{(h,g)},\ \
\ka\in \BZ^{(h,h)}\, \}.$$ \indent The aim of this paper is to
find a fundamental domain for $\G_{g,h}\ba \BH_{g,h}.$ This
article is organized as follows. In Section 2, we review the
Minkowski domain and the Siegel's fundamental domain $\CCF$
roughly. In Section 3, we find a fundamental domain for
$\G_{g,h}\ba \BH_{g,h}$ and present Riemannian metrics on the
fundamental domain invariant under the action (1.2) of the Jacobi
group $G^J.$ In Section 4, we investigate the spectral theory of
the Laplacian on the abelian variety $A_{\Om}$ associated to $\Om
\in \CCF$.

\end{section}
%
%
\begin{section}{{\bf Review on a Fundamental Domain $\CCF$ for $\G_g\ba \BH_g$}}
\setcounter{equation}{0} We let
$$\CP=\left\{\, Y\in\BR^{(g,g)}\,|\ Y=\,^tY>0\ \right\}$$
be an open cone in $\BR^N$ with $N=g(g+1)/2.$ The general linear
group $GL(g,\BR)$ acts on $\CP$ transitively by
\begin{equation}
g\circ Y:=gY\,^tg,\qquad g\in GL(g,\BR),\ Y\in \CP.\end{equation}
Thus $\CP$ is a symmetric space diffeomorphic to $GL(g,\BR)/O(g).$
For a matrix $A\in F^{(k,l)}$ and $B\in F^{(k,l)},$ we write
$A[B]=\,^tBAB$ and for a square matrix $A$, $\sigma(A)$ denotes
the trace of $A.$.
\vskip 0.10cm
\newcommand\Mg{{\mathcal M}_g}
\newcommand\Rg{{\mathcal R}_g}
The fundamental domain $\Rg$ for $GL(g,\BZ)\ba \CP$ which was
found by H. Minkowski\,[5] is defined as a subset of $\CP$
consisting of $Y=(y_{ij})\in \CP$ satisfying the following
conditions (M.1)-(M.2)\ (cf. [2] p.\,191 or [4] p.\,123): \vskip
0.1cm (M.1)\ \ \ $aY\,^ta\geq y_{kk}$\ \ for every
$a=(a_i)\in\BZ^g$ in which $a_k,\cdots,a_g$ are relatively prime
for $k=1,2,\cdots,g$. \vskip 0.1cm (M.2)\ \ \ \ $y_{k,k+1}\geq 0$
\ for $k=1,\cdots,g-1.$ \vskip 0.1cm We say that a point of $\Rg$
is {\it Minkowski reduced} or simply {\it M}-{\it reduced}. $\Rg$
has the following properties (R1)-(R6): \vskip 0.1cm (R1) \ For
any $Y\in\CP,$ there exist a matrix $A\in GL(g,\BZ)$ and $R\in\Rg$
such that $Y=R[A]$\,(cf. [2] p.\,191 or [4] p.\,139). That is,
$$GL(g,\BZ)\circ \Rg=\CP.$$
\indent (R2)\ \ $\Rg$ is a convex cone through the origin bounded
by a finite number of hyperplanes. $\Rg$ is closed in $\CP$
(cf.\,[4] p.\,139).

\vskip 0.1cm (R3) If $Y$ and $Y[A]$ lie in $\Rg$ for $A\in
GL(g,\BZ)$ with $A\neq \pm I_g,$ then $Y$ lies on the boundary
$\partial \Rg$ of $\Rg$. Moreover $\Rg\cap (\Rg [A])\neq
\emptyset$ for only finitely many $A\in GL(g,\BZ)$ (cf.\,[4]
p.\,139). \vskip 0.1cm (R4) If $Y=(y_{ij})$ is an element of
$\Rg$, then
$$y_{11}\leq y_{22}\leq \cdots \leq y_{gg}\quad \text{and}\quad
|y_{ij}|<{\frac 12}y_{ii}\quad \text{for}\ 1\leq i< j\leq g.$$
\indent We refer to [2] p.\,192 or [4] pp.\,123-124.  \vskip
0.1cm\noindent {\it Remark.} Grenier\,[1] found another
fundamental domain for $GL(g,\BZ)\ba \CP.$ \vskip 0.1cm
\newcommand\POB{ {{\partial}\over {\partial{\overline \Omega}}} }
\newcommand\PZB{ {{\partial}\over {\partial{\overline Z}}} }
\newcommand\PX{ {{\partial}\over{\partial X}} }
\newcommand\PY{ {{\partial}\over {\partial Y}} }
\newcommand\PU{ {{\partial}\over{\partial U}} }
\newcommand\PV{ {{\partial}\over{\partial V}} }
\newcommand\PO{ {{\partial}\over{\partial \Omega}} }
\newcommand\PZ{ {{\partial}\over{\partial Z}} }
\vskip 0.2cm For $Y=(y_{ij})\in \CP,$ we put
$$dY=(dy_{ij})\qquad\text{and}\qquad \PY\,=\,\left(\,
{ {1+\delta_{ij}}\over 2}\, { {\partial}\over {\partial y_{ij} } }
\,\right).$$ Then we can see easily that
\begin{equation}
ds^2=\s ( (Y^{-1}dY)^2)\end{equation} is a $GL(g,\BR)$-invariant
Riemannian metric on $\CP$ and its Laplacian is given by
$$\Delta=\s \left( \left( Y\PY\right)^2\right).$$
We also can see that
$$d\mu_g(Y)=(\det Y)^{-{ {g+1}\over2 } }\prod_{i\leq j}dy_{ij}$$
is a $GL(g,\BR)$-invariant volume element on $\CP$. The metric
$ds^2$ on $\CP$ induces the metric $ds_{\mathcal R}^2$ on $\Rg.$
Minkowski [5] calculated the volume of $\Rg$ for the volume
element $[dY]:=\prod_{i\leq j}dy_{ij}$ explicitly. Later Siegel
[7],\,[9] computed the volume of $\Rg$ for the volume element
$[dY] $ by a simple analytic method and generalized this case to
the case of any algebraic number field. \vskip 0.1cm Siegel\,[8]
determined a fundamental domain $\CCF$ for $\G_g\ba \BH_g.$ We say
that $\Om=X+iY\in \BH_g$ with $X,\,Y$ real is {\it Siegel reduced}
or {\it S}-{\it reduced} if it has the following three properties:
\vskip 0.1cm (S.1)\ \ \ $\det (\text{Im}\,(\g\cdot\Om))\leq \det
(\text{Im}\,(\Om))\qquad\text{for\ all}\ \g\in\G_g$; \vskip 0.1cm
(S.2)\ \ $Y=\text{Im}\,\Om$ is M-reduced, that is, $Y\in \Rg\,;$
\vskip 0.1cm (S.3) \ \ $|x_{ij}|\leq {\frac 12}\quad \text{for}\
1\leq i,j\leq g,\ \text{where}\ X=(x_{ij}).$ \vskip 0.1cm $\CCF$
is defined as the set of all Siegel reduced points in $\BH_g.$
Using the highest point method, Siegel proved the following
(F1)-(F3)\,(cf. [2] pp.\,194-197 or [4] p.\,169): \vskip 0.1cm
(F1)\ \ \ $\G_g\cdot \CCF=\BH_g,$ i.e.,
$\BH_g=\cup_{\g\in\G_g}\g\cdot \CCF.$ \vskip 0.1cm (F2)\ \ $\CCF$
is closed in $\BH_g.$ \vskip 0.1cm (F3)\ \ $\CCF$ is connected and
the boundary of $\CCF$ consists of a finite number of hyperplanes.
\vskip 0.21cm For $\Om=(\omega_{ij})\in\BH_g,$ we write $\Om=X+iY$
with $X=(x_{ij}),\ Y=(y_{ij})$ real and $d\Om=(d\om_{ij})$. We
also put
$$\PO=\,\left(\,
{ {1+\delta_{ij}}\over 2}\, { {\partial}\over {\partial \om_{ij} }
} \,\right) \qquad\text{and}\qquad \POB=\,\left(\, {
{1+\delta_{ij}}\over 2}\, { {\partial}\over {\partial {\overline
{\om}}_{ij} } } \,\right).$$ Then
\begin{equation}
ds_*^2=\s (Y^{-1}d\Om\, Y^{-1}d{\overline\Om})\end{equation} is a
$Sp(g,\BR)$-invariant K{\"a}hler metric on $\BH_g$ (cf.\,[8]) and
H. Maass [3] proved that its Laplacian is given by
\begin{equation}
\Delta_*=\,4\,\s \left( Y\,\,{
}^t\!\left(Y\POB\right)\PO\right).\end{equation} And
\begin{equation}
dv_g(\Om)=(\det Y)^{-(g+1)}\prod_{1\leq i\leq j\leq g}dx_{ij}\,
\prod_{1\leq i\leq j\leq g}dy_{ij}\end{equation} is a
$Sp(g,\BR)$-invariant volume element on
$\BH_g$\,(cf.\,[10],\,p.\,130). The metric $ds_*^2$ given by (2.3)
induces a metric $ds_{\mathcal F}^2$ on $\CCF.$ \vskip 0.1cm
Siegel\,[8] computed the volume of $\CCF$
\begin{equation}
\text{vol}\,(\CCF)=2\prod_{k=1}^g\pi^{-k}\G
(k)\zeta(2k),\end{equation} where $\G (s)$ denotes the Gamma
function and $\zeta (s)$ denotes the Riemann zeta function. For
instance,
$$\text{vol}\,({\mathcal F}_1)={{\pi}\over 3},\quad \text{vol}\,({\mathcal F}_2)={{\pi^3}\over {270}},
\quad \text{vol}\,({\mathcal F}_3)={{\pi^6}\over {127575}},\quad
\text{vol}\,({\Cal F}_4)={{\pi^{10}}\over {200930625}}.$$
\end{section}
%
%
\begin{section}{{\bf A Fundamental Domain for $\Gamma_{g,h}\backslash
\BH_{g,h}$}} \setcounter{equation}{0} Let $E_{kj}$ be the $h\times
g$ matrix with entry 1 where the $k$-th row and the $j$-th colume
meet, and all other entries 0. For an element $\Om\in \BH_g$, we
set for brevity
\begin{equation}
F_{kj}(\Om):=E_{kj}\Om,\qquad 1\leq k\leq h,\ 1\leq j\leq
n.\end{equation}
 \indent For each $\Om\in {\mathcal F}_g,$ we define a
subset $P_{\Om}$ of $\BC^{(h,g)}$ by
\begin{equation*}
P_{\Om}=\left\{ \,\sum_{k=1}^h\sum_{j=1}^g \la_{kj}E_{kj}+
\sum_{k=1}^h\sum_{j=1}^g \mu_{kj}F_{kj}(\Om)\,\Big|\ 0\leq
\la_{kj},\mu_{kj}\leq 1\,\right\}. \end{equation*} \indent For
each $\Om\in \CCF,$ we define the subset $D_{\Om}$ of $\BH_{g,h}$
by
\begin{equation*} D_{\Om}:=\left\{\,(\Om,Z)\in\BH_{g,h}\,\vert\ Z\in
P_{\Om}\,\right\}.\end{equation*}
\newcommand\Fgh{{\mathcal F}_{g,h}}
We define
\begin{equation*} \Fgh:=\cup_{\Om\in\CCF}D_{\Omega}.\end{equation*}
\begin{theorem}
$\Fgh$ is a fundamental domain for $\G_{g,h}\ba \BH_{g,h}.$
\end{theorem}
\begin{proof} Let $({\tilde{\Om}},{\tilde{Z}})$ be an arbitrary
element of $\BH_{g,h}.$ We must find an element $(\Om,Z)$ of
$\Fgh$ and an element $\g^J=(\g,(\la,\mu;\ka))\in\G_{g,h}$ with
$\g\in\G_g$ such that $\g^J\cdot
(\Om,Z)=({\tilde{\Om}},{\tilde{Z}}).$ Since $\CCF$ is a
fundamental domain for $\G_g\ba \BH_g,$ there exists an element
$\g$ of $\G_g$ and an element $\Om$ of $\CCF$ such that
$\g\cdot\Om={\tilde {\Om}}.$ Here $\Om$ is unique up to the
boundary of $\CCF$. \vskip 0.1cm We write
$$\g=\begin{pmatrix} A & B\\ C & D \end{pmatrix} \in \G_g.$$
It is easy to see that we can find $\la,\mu\in \BZ^{(h,g)}$ and
$Z\in P_{\Om}$ satisfying the equation
$$Z+\la \Om +\mu={\tilde Z}(C\Om+D).$$
If we take $\g^J=(\g,(\la,\mu;0))\in \G_{g,h},$ we see that
$\g^J\cdot (\Om,Z)=({\tilde{\Om}},{\tilde{Z}}).$ Therefore we
obtain
$$\BH_{g,h}=\cup_{\g^J\in \G_{g,h}}\g^J\cdot\Fgh.$$
Let $(\Om,Z)$ and $\g^J\cdot (\Om,Z)$ be two elements of $\Fgh$
with $\g^J=(\g,(\la,\mu;\ka))\in \G_{g,h}.$ Then both $\Om$ and
$\g\cdot\Om$ lie in $\CCF$. Therefore both of them either lie in
the boundary of $\CCF$ or $\g=\pm I_{2g}.$ In the case that both
$\Om$ and $\g\cdot\Om$ lie in the boundary of $\CCF$, both
$(\Om,Z)$ and $\g^J\cdot (\Om,Z)$ lie in the boundary of $\Fgh$.
If $\g=\pm I_{2g},$ we have
\begin{equation}
Z\in P_{\Om}\quad \text{and}\quad \pm (Z+\la \Om+\mu)\in
P_{\Om},\quad \la,\mu\in \BZ^{(h,g)}.\end{equation} From the
definition of $P_{\Om}$ and (3.2), we see that either $\la=\mu=0,\
\g\neq -I_{2g}$ or both $Z$ and $\pm (Z+\la \Om+\mu)$ lie on the
boundary of the parallelepiped $P_{\Om}$. Hence either
both$(\Om,Z)$ and $\g^J\cdot (\Om,Z)$ lie in the boundary of
$\Fgh$ or $\g^J=(I_{2g},(0,0;\ka))\in\G_{g,h}$. Consequently
$\Fgh$ is a fundamental domain for $\G_{g,h}\ba \BH_{g,h}.$
\end{proof} For a coordinate $(\Om,Z)\in\BH_{g,h}$ with
$\Om=(\om_{\mu\nu})\in {\mathbb H}_g$ and $Z=(z_{kl})\in \Chg,$ we
put
\begin{align*}
\Om\,=&\,X\,+\,iY,\quad\ \ X\,=\,(x_{\mu\nu}),\quad\ \
Y\,=\,(y_{\mu\nu})
\ \ \text{real},\\
Z\,=&U\,+\,iV,\quad\ \ U\,=\,(u_{kl}),\quad\ \ V\,=\,(v_{kl})\ \
\text{real},\\
d\Om\,=&\,(d\om_{\mu\nu}),\quad\ \ dX\,=\,(dx_{\mu\nu}),\quad\ \
dY\,=\,(dy_{\mu\nu}),\\
dZ\,=&\,(dz_{kl}),\quad\ \ dU\,=\,(du_{kl}),\quad\ \
dV\,=\,(dv_{kl}),\\
d{\overline{\Om}}=&\,(d{\overline{\om}}_{\mu\nu}),\quad
d{\overline Z}=(d{\bar z}_{kl}),
\end{align*}
\newcommand\POB{ {{\partial}\over {\partial{\overline \Omega}}} }
\newcommand\PZB{ {{\partial}\over {\partial{\overline Z}}} }
\newcommand\PX{ {{\partial}\over{\partial X}} }
\newcommand\PY{ {{\partial}\over {\partial Y}} }
\newcommand\PU{ {{\partial}\over{\partial U}} }
\newcommand\PV{ {{\partial}\over{\partial V}} }
\newcommand\PO{ {{\partial}\over{\partial \Omega}} }
\newcommand\PZ{ {{\partial}\over{\partial Z}} }
$$ {{\partial}\over{\partial \Omega}}\,=\,\left(\, { {1+\delta_{\mu\nu}} \over 2}\, {
{\partial}\over {\partial \om_{\mu\nu}} } \,\right),\quad
\POB\,=\,\left(\, { {1+\delta_{\mu\nu}}\over 2} \, {
{\partial}\over {\partial {\overline \om}_{\mu\nu} }  }
\,\right),$$
$$\PZ=\begin{pmatrix} { {\partial}\over{\partial z_{11}} } & \hdots &
{ {\partial}\over{\partial z_{h1}} }\\
\vdots&\ddots&\vdots\\
{ {\partial}\over{\partial z_{1g}} }&\hdots &{ {\partial}\over
{\partial z_{hg}} } \end{pmatrix},\quad \PZB=\begin{pmatrix} {
{\partial}\over{\partial {\overline z}_{11} }   }&
\hdots&{ {\partial}\over{\partial {\overline z}_{h1} }  }\\
\vdots&\ddots&\vdots\\
{ {\partial}\over{\partial{\overline z}_{1g} }  }&\hdots & {
{\partial}\over{\partial{\overline z}_{hg} }  }
\end{pmatrix}.$$

{\it Remark.} The following metric \begin{align*}
ds_{g,h}^2=&\,\s\left(Y^{-1}d\Om\,Y^{-1}d{\overline \Om}\right)\,+
\,\s\left(Y^{-1}\,^tV\,V\,Y^{-1}d\Om\,Y^{-1}
d{\overline{\Om}}\right) \notag\\
&\ \ \ \ +\,\s\left(Y^{-1}\,^t(dZ)\,d{\overline Z}\right)\\
&\ \ -\s\left(\,V\,Y^{-1}d\Om\,Y^{-1}\,^t( d{\overline{\Om}} )\,
+\,V\,Y^{-1} d{\overline {\Om}} \, Y^{-1}\,^t(dZ)\,\right)\notag
\end{align*}
is a K{\"a}hler metric on $\BH_{g,h}$ which is invariant under the
action (1.2) of the Jacobi group $G^J$. Its Laplacian is given by
\begin{align*}
\Delta_{g,h}\,=\,& 4\,\s\left(\,Y\,\,
^t\!\left(Y\POB\right)\PO\,\right)\,+\,
4\,\s\left(\, Y\,\PZ\,\,{}^t\!\left( \PZB\right)\,\right) \nonumber \\
&\ \ \ \ +\,4\,\s\left(\,VY^{-1}\,^tV\,\,^t\!\left(Y\PZB\right)\,\PZ\,\right)\\
&\ \
+\,4\,\s\left(V\,\,^t\!\left(Y\POB\right)\PZ\,\right)+\,4\s\left(\,^tV\,\,^t\!\left(Y\PZB\right)\PO\,\right).\nonumber
\end{align*}

The following differential form
$$dv_{g,h}=\,\left(\,\text{det}\,Y\,\right)^{-(g+h+1)}[dX]\w [dY]\w
[dU]\w [dV]$$ is a $G^J$-invariant volume element on $\BH_{g,h}$,
where
$$[dX]=\w_{\mu\leq\nu}dx_{\mu\nu},\quad [dY]=\w_{\mu\leq\nu}
dy_{\mu\nu},\quad [dU]=\w_{k,l}du_{kl}\quad \text{and} \quad
[dV]=\w_{k,l}dv_{kl}.$$ The point is that the invariant metric
$ds_{g,h}^2$ and its Laplacian are beautifully expressed in terms
of the {\it trace} form. The proof of the above facts can be found
in [15].

\newcommand\Imm{\text{Im}\,}
\newcommand\MCM{\mathcal M}
\end{section}
%
%
\newcommand\Fgh{{\mathcal F}_{g,h}}
\begin{section}{{\bf Spectral Decomposition of $L^2(A_{\Om}$)}}
\setcounter{equation}{0}  \vskip 0.2cm We fix two positive
integers $g$ and $h$ throughout this section. \vskip 0.1cm For an
element $\Om\in \BH_g,$ we set
\begin{equation*}
L_{\Om}:=\BZ^{(h,g)}+\BZ^{(h,g)}\Om\end{equation*} We use the
notation (3.1). It follows from the positivity of $\text{Im}\,\Om$
that the elements $E_{kj},\,F_{kj}(\Om)\\ (1\leq k\leq h,\ 1\leq
j\leq g)$ of $L_{\Om}$ are linearly independent over $\BR$.
Therefore $L_{\Om}$ is a lattice in $\BC^{(h,g)}$ and the set
$\left\{\,E_{kj},\,F_{kj}(\Om)\,|\ 1\leq k\leq h,\ 1\leq j\leq g\,
\right\}$ forms an integral basis of $L_{\Om}$. We see easily that
if $\Om$ is an element of $\BH_g$, the period matrix
$\Om_*:=(I_g,\Om)$ satisfies the Riemann conditions (RC.1) and
(RC.2)\,: \vskip 0.1cm (RC.1) \ \ \ $\Om_*J_g\,^t\Om_*=0\,$;
\vskip 0.1cm (RC.2) \ \ \ $-{1 \over
{i}}\Om_*J_g\,^t{\overline{\Om}}_*
>0$.

\vskip 0.2cm \noindent Thus the complex torus
$A_{\Om}:=\BC^{(h,g)}/L_{\Omega}$ is an abelian variety. For more
details on $A_{\Om}$, we refer to [2] and [6]. \vskip 0.2cm It
might be interesting to investigate the spectral theory of the
Laplacian $\Delta_{g,h}$ on a fundamental domain $\Fgh$. But this
work is very complicated and difficult at this moment. It may be
that the first step is to develop the spectral theory of the
Laplacian $\Delta_{\Omega}$ on the abelian variety $A_{\Omega}.$
The second step will be to study the spectral theory of the
Laplacian $\Delta_*$\,(see (2.4)) on the moduli space
$\Gamma_g\backslash \BH_g$ of principally polarized abelian
varieties of dimension $g$. The final step would be to combine the
above steps and more works to develop the spectral theory of the
Lapalcian $\Delta_{g,h}$ on $\Fgh.$ In this section, we deal only
with the spectral theory $\Delta_{\Omega}$ on $L^2(A_{\Omega}).$

 \vskip 0.1cm We fix
an element $\Om=X+iY$ of $\BH_g$ with $X=\text{Re}\,\Om$ and
$Y=\text{Im}\, \Om.$ For a pair $(A,B)$ with $A,B\in\BZ^{(h,g)},$
we define the function $E_{\Om;A,B}:\Chg\lrt \BC$ by
\begin{equation*}
E_{\Om;A,B}(Z)=e^{2\pi i\left( \s\,(\,^tAU\,)+\,\s\,
((B-AX)Y^{-1}\,^tV)\right)},\end{equation*} where $Z=U+iV$ is a
variable in $\Chg$ with real $U,V$. \vskip 0.1cm\noindent
\begin{lemma} For any $A,B\in \BZ^{(h,g)},$ the function
$E_{\Om;A,B}$ satisfies the following functional equation
\begin{equation*}
E_{\Om;A,B}(Z+\la \Om+\mu)=E_{\Om;A,B}(Z),\quad
Z\in\Chg\end{equation*} for all $\la,\mu\in\BZ^{(h,g)}.$ Thus
$E_{\Om;A,B}$ can be regarded as a function on $A_{\Om}.$ \vskip
0.1cm \end{lemma}
\begin{proof}
We write $\Om=X+iY$ with real $X,Y.$ For any
$\la,\mu\in\BZ^{(h,g)},$ we have
\begin{align*}
E_{\Om;A,B}(Z+\la\Om+\mu)&=E_{\Om;A,B}((U+\la X+\mu)+i(V+\la Y))\\
&=e^{ 2\pi i \left\{\,\s\,(\,^t\!A(U+\la
X+\mu))+\,\s\,((B-AX)Y^{-1}\,^t\!(V+\la Y))\,\right\} }\\
&=e^{ 2\pi i \left\{\,\s\,(\,^t\!AU+\,^t\!A\la X+\,^t\!A\mu)
+\,\s\,((B-AX)Y^{-1}\,^tV+B\,^t\la-AX\,^t\la) \right\} }\\
&=e^{2\pi i \left\{\,\s\,(\,^t\!AU)\,+\,\s\,((B-AX)Y^{-1}\,^tV)\right\} }\\
&=E_{\Om;A,B}(Z).\end{align*} Here we used the fact that
$^t\!A\mu$ and $B\,^t\la$ are integral. \end{proof}
\newcommand\AO{A_{\Omega}}
\newcommand\Imm{\text{Im}}
We use the notations in Section 3.
\begin{lemma}
The metric
$$ds_{\Om}^2=\s\left((\textrm{Im}\,\Om)^{-1}\,\,^t(dZ)\,d{\overline Z})\,\right)$$
is a K{\"a}hler metric on $A_{\Om}$ invariant under the action
(1.2) of $\G^J=Sp(g,\BZ)\ltimes H_{\BZ}^{(h,g)}$ on $(\Om,Z)$ with
$\Om$ fixed. Its Laplacian $\Delta_{\Om}$ of $ds_{\Om}^2$ is given
by
\begin{equation*}
\Delta_{\Om}=\,\s\left( (\Imm\,\Omega)\,{ {\partial}\over
{\partial {Z}} }\,^t\!\left(  {{\partial}\over {\partial
{\overline Z}}} \right)\,
 \right). \end{equation*}
\end{lemma}
\begin{proof} Let ${\tilde \gamma}=(\gamma,(\la,\mu;\kappa))\in
\Gamma^J$ with $\gamma=\begin{pmatrix} A & B\\ C & D
\end{pmatrix}\in Sp(g,\BZ)$ and $({\tilde \Omega},{\tilde
Z})={\tilde \gamma}\cdot (\Omega,Z)$ with $\Omega\in {\mathbb
H}_g$ fixed. Then according to [4],\,p.\,33,
$$\Imm\,\g \cdot\Omega=\,^t(C{\overline
\Om}+D)^{-1}\,\Imm\,\Om\,(C\Om+D)^{-1}$$ and by (1.2),
$$d{\tilde Z}=dZ\,(C\Om+D)^{-1}.$$
Therefore
\begin{eqnarray*}
& & (\Imm\,{\tilde\Om})^{-1}\,^t(d{\tilde Z})\,d{\overline{\tilde Z}}  \\
&=&(C{\overline\Om}+D)\, (\Imm\,\Om)^{-1}\,^t(C
{\Om}+D)\,^t(C{ \Om}+D)^{-1}\,^t(d{ Z})\,d{\overline Z}\,(C{\overline \Om}+D)^{-1} \\
&=& (C{\overline\Om}+D)\,(\Imm\,\Om)^{-1}\,^t(dZ)\, d{\overline
Z}\,(C{\overline \Om}+D)^{-1}.
\end{eqnarray*}
The metric $ds_{iI_g}=\s (dZ\,^t(d{\overline Z}))$ at $Z=0$ is
positive definite. Since $G^J$ acts on ${\mathbb H}_{g,h}$
transitively, $ds^2_{\Om}$ is a Riemannian metric for any $\Om\in
{\mathbb H}_g.$ We note that the differential operator
$\Delta_{\Om}$ is invariant under the action of $\Gamma^J.$ In
fact, according to (1.2),
$${ {\partial}\over {\partial {\tilde Z}} }=(C\Om +D)\,{
{\partial}\over {\partial Z} }.$$ Hence if $f$ is a differentiable
function on $A_{\Om}$, then
\begin{eqnarray*}
& & \Imm\,{\tilde \Omega}\, { {\partial}\over {\partial {\tilde
Z}} }
\,^t\!\left( { {\partial f}\over {\partial \overline{\tilde Z} }  } \right) \\
&=&\,^t(C{\overline\Om}+D)^{-1}\,(\Imm\,\Om)\,(C{\Om}+D)^{-1} (C{
\Om}+D)\, { {\partial}\over {\partial Z} }\,^t\!\left(
(C{\overline\Om}+D){ {\partial f}\over {\partial \overline Z} }
\right)\\
&=&\,^t(C{\overline\Om}+D)^{-1} \, \Imm\,\Omega\,{ {\partial}\over
{\partial Z} }\,^t\!\left( { {\partial f}\over {\partial \overline
Z} } \right)(C{\overline\Om}+D).
\end{eqnarray*}
Therefore
\begin{equation*}
\s\left(\Imm\,{\tilde \Omega}\, { {\partial}\over {\partial
{\tilde Z}} } \,^t\!\left( { {\partial }\over {\partial
\overline{\tilde Z} }  } \right) \right)=\,\sigma\left(\,
\Imm\,\Omega\,{ {\partial}\over {\partial Z} }\,^t\!\left( {
{\partial f}\over {\partial \overline Z} } \right) \right).
\end{equation*}

 By the induction on $h$, we can compute the Laplacian
$\Delta_{\Om}.$

\end{proof}
  \vskip 0.1cm We
let $L^2(\AO)$ be the space of all functions $f:\AO\lrt\BC$ such
that
$$||f||_{\Om}:=\int_{\AO}|f(Z)|^2dv_{\Om},$$
where $dv_{\Om}$ is the volume element on $\AO$ normalized so that
$\int_{\AO}dv_{\Om}=1.$ The inner product $(\,\,,\,\,)_{\Om}$ on
the Hilbert space $L^2(\AO)$ is given by
\begin{equation}
(f,g)_{\Om}:=\int_{\AO}f(Z)\,{\overline{g(Z)} }\,dv_{\Om},\quad
f,g\in L^2(\AO).\end{equation}
\begin{theorem}
The set $\left\{\,E_{\Om;A,B}\,|\ A,B\in\BZ^{(h,g)}\,\right\}$ is
a complete orthonormal basis for $L^2(\AO)$. Moreover we have the
following spectral decomposition of $\Delta_{\Om}$:
$$L^2(\AO)=\oplus_{A,B\in \BZ^{(h,g)}}\BC\cdot E_{\Om;A,B}.$$
\end{theorem}
\begin{proof} Let
\begin{equation*}
T=\BC^{(h,g)}/(\BZ^{(h,g)}\times \BZ^{(h,g)})=(\BR^{(h,g)}\times
\BR^{(h,g)})/(\BZ^{(h,g)}\times \BZ^{(h,g)})
\end{equation*}
be the torus of real dimension $2hg$. The Hilbert space $L^2(T)$
is isomorphic to the $2hg$ tensor product of $L^2(\BR/\BZ)$, where
$\BR/\BZ$ is the one-dimensional real torus. Since
$L^2(\BR/\BZ)=\oplus_{n\in\BZ}\BC\cdot e^{2\pi inx},$ the Hilbert
space $L^2(T)$ is
\begin{equation*}
L^2(T)=\oplus_{A,B\in \BZ^{(h,g)}}\BC\cdot E_{A,B}(W),
\end{equation*}
where $W=P+iQ,\ P,Q\in \BR^{(h,g)}$ and
\begin{equation*}
E_{A,B}(W):=e^{2\pi i\,\sigma(\,^t\!AP+\,^tBQ)},\quad
A,B\in\BZ^{(h,g)}.
\end{equation*}
The inner product on $L^2(T)$ is defined by
\begin{equation}
(f,g):=\int_0^1\cdots \int_0^1 f(W)\,{\overline{g(W)}
}\,dp_{11}\cdots dp_{hg}dq_{11}\cdots dq_{hg},\quad f,g\in
L^2(T),\end{equation}
where $W=P+iQ\in T,\ P=(p_{kl})$ and
$Q=(q_{kl}).$ Then we see that the set $\left\{ E_{A,B}(W)\,|\
A,B\in\BZ^{(h,g)}\,\right\}$ is a complete orthonormal basis for
$L^2(T)$, and each $E_{A,B}(W)$ is an eigenfunction of the
standard Laplacian
\begin{equation*}
\Delta_T=\sum_{k=}^h\sum_{l=1}^g \left( {{\partial^2}\over
{\partial p_{kl}^2} }+{{\partial^2}\over {\partial q_{kl}^2}
}\right).
\end{equation*}
We define the mapping $\Phi_{\Omega}:T\lrt A_{\Omega}$ by
\begin{equation}
\Phi_{\Omega}(P+iQ)=(P+QX)+i\,QY, \quad P+i\,Q\in T,\ P,Q\in
\BR^{(h,g)}.
\end{equation}
This is well defined. We can see that $\Phi_{\Omega}$ is a
diffeomorphism and that the inverse $\Phi_{\Omega}^{-1}$ of
$\Phi_{\Omega}$ is given by
\begin{equation}
\Phi_{\Omega}^{-1}(U+iV)=(U-VY^{-1}X)\,+\,i\,VY^{-1}, \quad
U+iV\in A_{\Omega},\ U,V\in \BR^{(h,g)}.
\end{equation}
Using (4.4), we can show that for $A,B\in\BZ^{(h,g)}$, the
function $E_{A,B}(W)$ on $T$ is transformed to the function
$E_{\Omega;A,B}$ on $A_{\Omega}$ via the diffeomorphism
$\Phi_{\Omega}$. Using (4.2) and the diffeomorphism
$\Phi_{\Omega}$, we can choose a normalized volume element
$dv_{\Omega}$ on $A_{\Omega}$ and then we get the inner product on
$L^2(A_{\Omega})$ defined by (4.1). This completes the proof.

\end{proof}

\end{section}

\vskip 0.1cm

\vspace{0.5cm}


\end{document}